\documentclass[a4paper,11pt]{article}

\usepackage{textcomp}
\usepackage{amsfonts}
\usepackage[dvips]{graphicx}
\usepackage{amsmath}

\author{Yi SONG \& Stephen P. BANKS\\Department of Automatic Control and Systems Engineering,
\\University of Sheffield, Mappin Street,\\Sheffield S1 3JD.
\\e-mail: s.banks@sheffield.ac.uk}

\title{\textbf{Dynamical Systems On Three Manifolds}\\\textbf{Part II:
3-Manifolds, \emph{Heegaard Splittings} and Three-Dimensional
Systems}}
\begin{document}

\maketitle

\newtheorem{theorem}{Theorem}[section]
\newtheorem{definition}{Definition}[section]
\newtheorem{proposition}{Proposition}[section]
\newtheorem{corollary}{Corollary}[section]
\newtheorem{lemma}{Lemma}[section]
\newtheorem{example}{Example}[section]

\begin{abstract}
The global behaviour of nonlinear systems is extremely important in
control and systems theory since the usual local theories will only
give information about a system in some neighbourhood of an
operating point. Away from that point, the system may have totally
different behaviour and so the theory developed for the local system
will be useless for the global one.

In this paper we shall consider the analytical and topological
structure of systems on 2- and 3- manifolds and show that it is
possible to obtain systems with 'arbitrarily strange' behaviour,
i.e., arbitrary numbers of chaotic regimes which are knotted and
linked in arbitrary ways. We shall do this by considering
\emph{Heegaard Splittings} of these manifolds and the resulting
systems defined on the boundaries.
\emph{C-homeomorphisms}.
\end{abstract}

\section{Introduction}
In a recent paper ([Banks \& Song, 2006]), we have shown how to
define general (analytic) systems on 2-manifolds by using the theory
of \emph{automorphic functions}. The importance of this theory to
dynamical systems is that, globally, they are defined not on `flat'
\emph{Euclidean} spaces, but on manifolds. In fact, it was shown in
([Banks \& Song, 2006]) that the simple pendulum `sits' naturally on
a \emph{Klein bottle}. In this paper, we consider the case of
three-dimensional systems and derive some results on the nature of
three-dimensional dynamical systems and the 3-manifolds on which
they `live'.

The main difficulty compared with the 2-manifold case is that
3-manifold topology is much more complex. Indeed, there is no
procedure for finding a complete set of topological invariants for a
three manifold although a great many invariants have been found,
surprisingly from quantum group theory ([Ohtsuki, 2001]). There we
shall extend our 2-manifold theory coupled with \emph{Heegaard
Splittings} and \emph{Connected Sums} to approach a theory of
3-dimensional dynamical systems.

\section{Three Manifolds and \emph{Heegaard Splittings}}
We shall consider, in this paper, dynamical systems defined on
3-manifolds. A 3-manifold \emph{M} is a separable metric space such
that each point $x \in \emph{M}$ has an open neighbourhood, which is
homeomorphic to $\mathbb{R}^3$ or $\mathbb{R}_{+}^3=\{x
\in{\mathbb{R}^3}: \; x_3\ge{0}\}$, we can assume all the
3-manifolds we consider here are differentiable (or
p.l.\footnote{piecewise linear}) manifolds since any 3-manifold has
a unique p.l. or differentiable structure (see [Hempel, 1976]).
Points in \emph{M} which look locally like $\mathbb{R}^3$ are called
boundary points. The set of all boundary points is denoted by
$\partial\emph{M}$. Note that $\partial \partial \emph{M} =
\emptyset$. A manifold which is compact and for which $\partial
\emph{M} = \emptyset$ is called \emph{closed}.

\begin{definition}
A Heegaard Splitting of a closed connected 3-manifold M is a pair
$(C_1, C_2)$ of cubes with handles such that
\[M=C_1\cup C_2\]
and
\[C_1 \cap C_2 =\partial C_1 = \partial C_2.\]
\end{definition}
The following results are well known (see, e.g. [Hempel, 1976]):
\begin{theorem}
Every closed, connected 3-manifold has a Heegaard Splitting.
\end{theorem}

\begin{proposition} \label{identifying klein bottle}
There is exactly one nonorientable 3-manifold with a genus one
Heegaard Splitting, the nonorientable 2-sphere bundle over $S^1$,
\emph{i.e.}, the trivial gluing of two solid Klein bottle.
\end{proposition}

Let $(C_1,C_2)$ be a \emph{Heegaard Splitting} of a 3-manifold
\emph{M}. A \emph{Heegaard diagram}, $(C_1; \partial D_1, \cdots,
\partial D_n)$, for the splitting $(C_1,C_2)$ consists of  a set
$\{D_1,\cdots, D_n\}$ of pairwise disjoint, properly embedded,
2-cells in $C_2$ which cut it into a 3-cell. We can regard \emph{M}
as being obtained from $C_1$ and $C_2$ by choosing a homeomorphism
of $\partial C_1$ onto $\partial C_2$ which maps a standard set of
longitudinal or meridian curves on $\partial C_1$ to $\{\partial
D_1, \cdots, \partial D_n\}$ situated on $\partial C_2$ (and
extending this homeomorphism throughout $C_1$ and $C_2$). Lickorish
([Lickorish, 1962]) shows that such a surface homeomorphism can be
generated (up to isotopy) by a sequence of \emph{C-homeomorphisms},
i.e., homeomorphisms of the following form:

Take a nontrivial cycle \emph{l} on the surface \emph{S}, cut along
\emph{l}, twist one side of the cycle through $2\pi$ and reconnect
the `two sides' of \emph{l}.

\begin{figure}[!hbp]
  \begin{center}
  \includegraphics[width=3in]{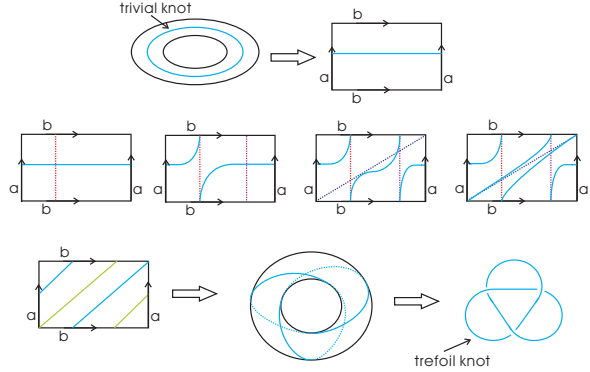}\\
  \caption{Generating a \emph{trefoil} knot from a trivial one via \emph{C-homeomorphisms}}\label{trefoil knot}
  \end{center}
\end{figure}

As an example, \textbf{Fig} \ref{trefoil knot} shows us how to get a
\emph{trefoil} knot from a trivial one in this way.

We shall use this result to perform surgery on our 2-dimensional
\emph{automorphic} systems and their extensions to obtain dynamical
systems on 3-manifolds in \S 3. In \S 4 we shall look for sufficient
conditions under which a nonlinear dynamical system on a 3-manifold
\emph{M} carries a \emph{Heegaard Splitting} which is compatible
with the dynamics in the sense that the \emph{Heegaard} surface is
invariant.

\section{Gluing Two Systems}
In this section we shall consider generating a three-dimensional
dynamical system by gluing together two systems defined on `cubes
with handles' along specified links. Modifying systems along links
to generate Pseudo-Anosov diffeomorphisms has been considered in
[Lozano, 1997]. Here we will apply the results of [Banks and Song,
2006] to generate an analytic (\emph{automorphic}) system on one
manifold and induce a twisted version on the other manifold by using
the so-called \emph{C-homeomorphisms} of [Lickorish, 1962].

\begin{figure}[!hbp]
  \begin{tabular}{cc}
  \includegraphics[width=2.2in]{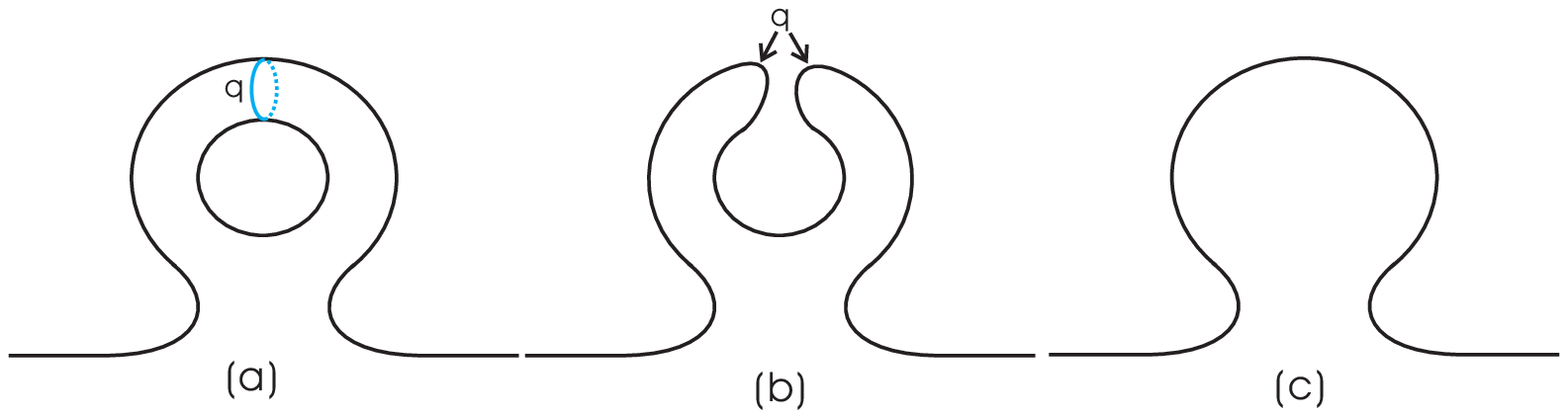} &
  \includegraphics[width=2.2in]{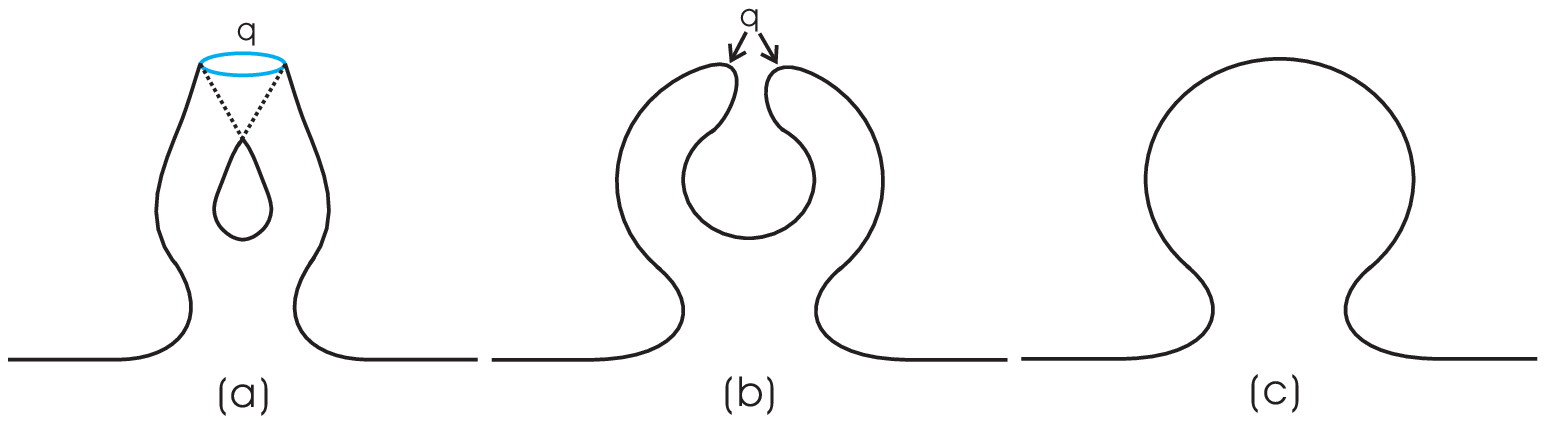}
  \end{tabular}
  \begin{center}
  \includegraphics[width=3in]{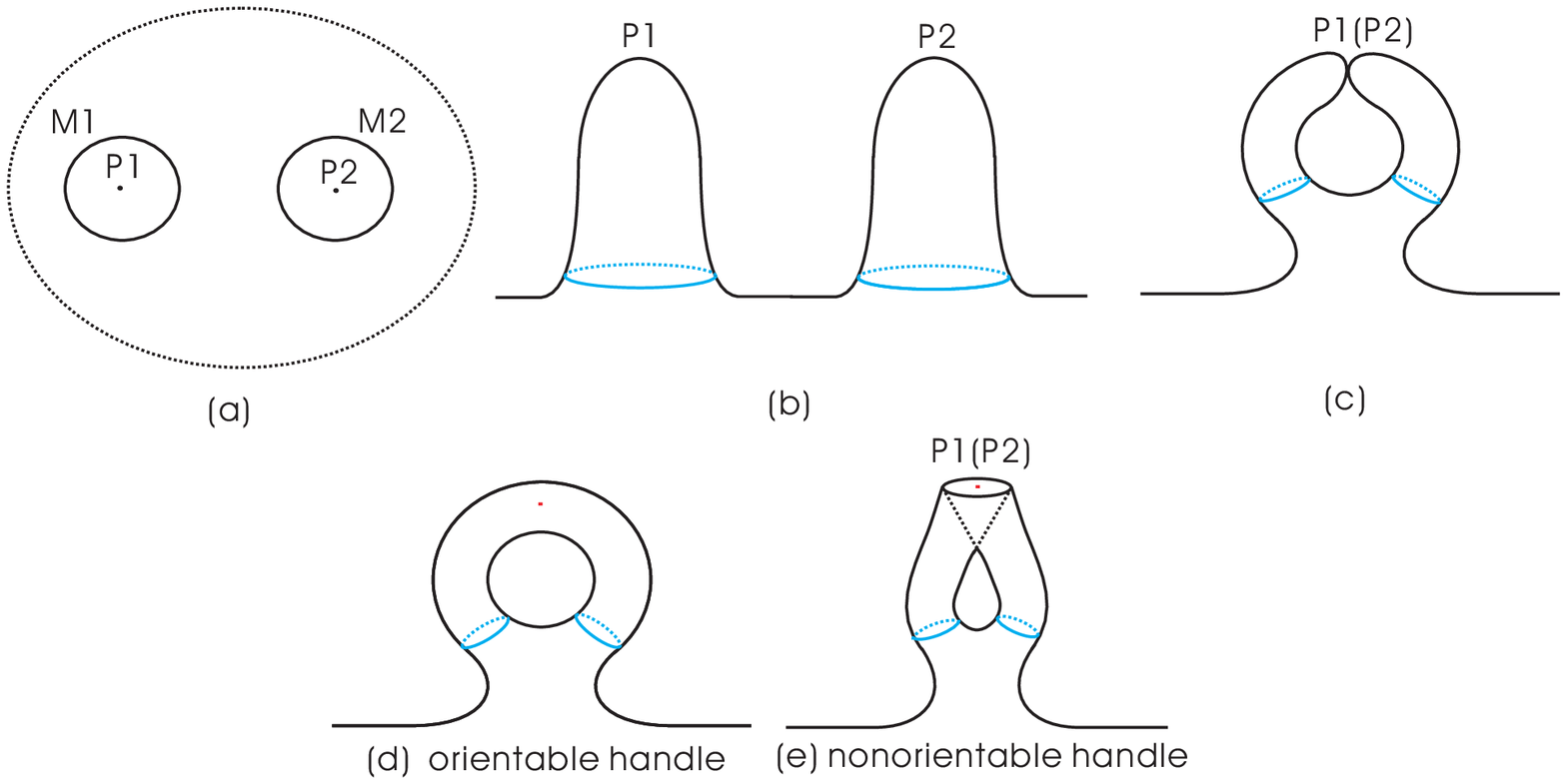}
  \end{center}
  \caption{Extending the surface dynamics throughout a solid handle}\label{handle}
\end{figure}

Suppose, therefore, that we wish to determine the analytic systems
defined on compact 3-manifolds which have an invariant surface
contained in the manifold. Let \emph{M} be a 3-manifold of that kind
with boundary \emph{S} which is a surface of \emph{genus g}. As
shown in [Banks and Song, 2006], a dynamical system on \emph{S} is
given by a \emph{generalized automorphic function F}, which
satisfies

\begin{equation} \label{auto}
F(Tz)= \frac{ad-bc}{(cz+d)^2}F(z),\qquad T \in \Gamma
\end{equation}
where $\Gamma$ is any \emph{Fuchsian group} and $T \in \Gamma$ is of
the form

\begin{equation}
T(z)= \frac{az+b}{cz+d}
\end{equation}
Any meromorphic function satisfying \textbf{Equation}(\ref{auto}) is
called an \emph{automorphic vector field} on \emph{S}. The neat
result shows that we can extend a meromorphic system defined on
\emph{S} as above to the whole of \emph{M} by adding a single
equilibrium point in $M/S$, plus one in each handle.
\begin{theorem}\label{extending dynamics}
Given a dynamical system on a surface S of genus g, we can extend it
to a dynamical system defined throughout the solid handle-body with
boundary S by adding a single equilibrium at the `centre' and one in
the interior of each handle.
\end{theorem}
\textbf{Proof.} Let $\{D_1, \cdots, D_g\}$ be a set of disjoint
properly embedded 2-cells in \emph{M} which cut \emph{M} into a ball
(3-cell) which do not contain any equilibriums on \emph{S}, and
shrink these 2-cells to points. We again obtain a 3-cell with $2g$
extra equilibrium points on the boundary. We may then regard this
3-cell as a standard ball with a spherical boundary. Now extend the
system defined on the surface into the whole 3-ball by simply
shrinking the surface dynamics to fit on a nested set of spheres
which fill out the 3-ball. thus the dynamics are foliated on
concentric spheres, and are identical on each sphere. The
singularity at the origin has index $2(1-g)$ by $Poincar\acute{e}$'s
theorem. To remove the equilibria inside the 3-ball apart from the
one at the origin, we add a normal vector field to the spheres which
is zero at the origin and the surface of the 3-ball and nonzero
elsewhere. Having defined an extension on the 3-ball we can return
to the original 3-manifold with a surface of \emph{genus g} by
gluing the appropriate points of the sphere and `blowing up' the
singularities there. This can clearly be done so that each resulting
handle has a single equilibrium in its interior. This process is
shown in \textbf{Fig} \ref{handle}. $\Box$

Now let us see some examples.
\begin{example}
\emph{A single pendulum is given by the following dynamical
equations}

\begin{eqnarray*}
\dot{\theta} &=& {\omega} \\
\dot{\omega} &=& -\frac{g}{l}\sin{\theta}
\end{eqnarray*}

\emph{\textbf{Fig} \ref{klein bottle}.(a) gives the dynamics in the
\emph{phase-plane}.}
\begin{figure}[!hbp]
\begin{center}
  \begin{tabular}{cc}
    \includegraphics[width=2.22in,height=2.2in]{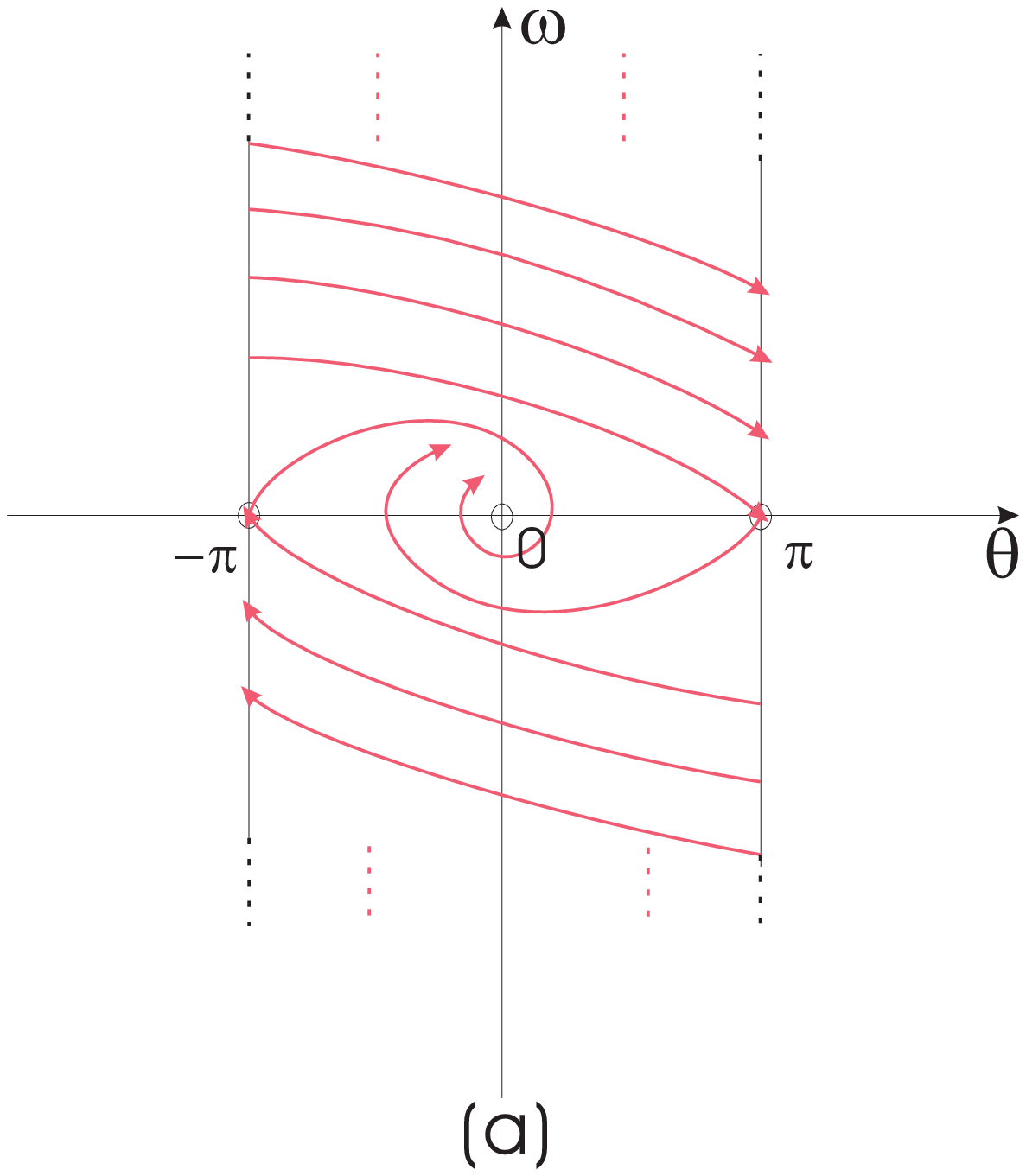} &
    \includegraphics[width=2.22in,height=2.2in]{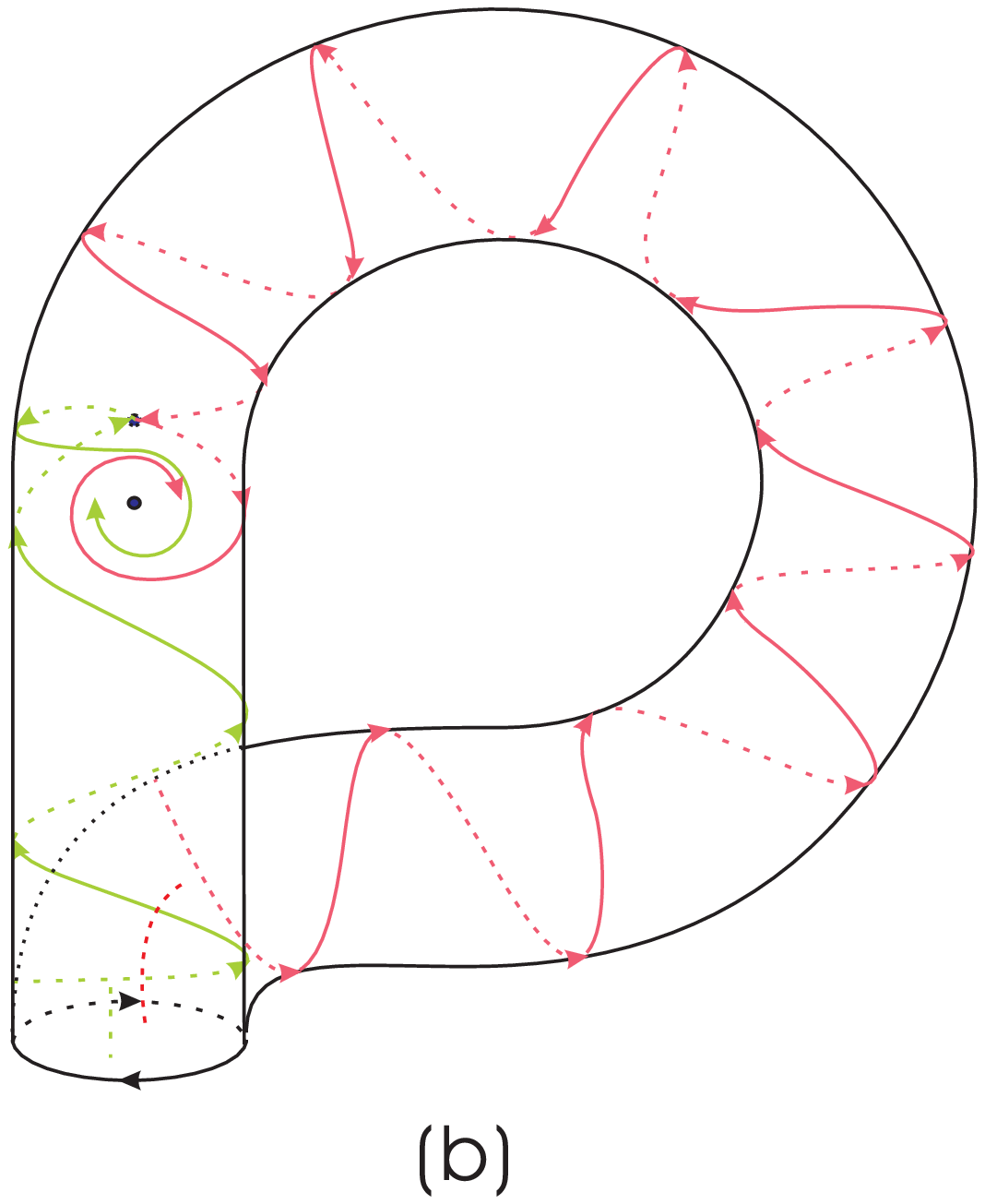}\\
    \includegraphics[width=2.22in,height=2.2in]{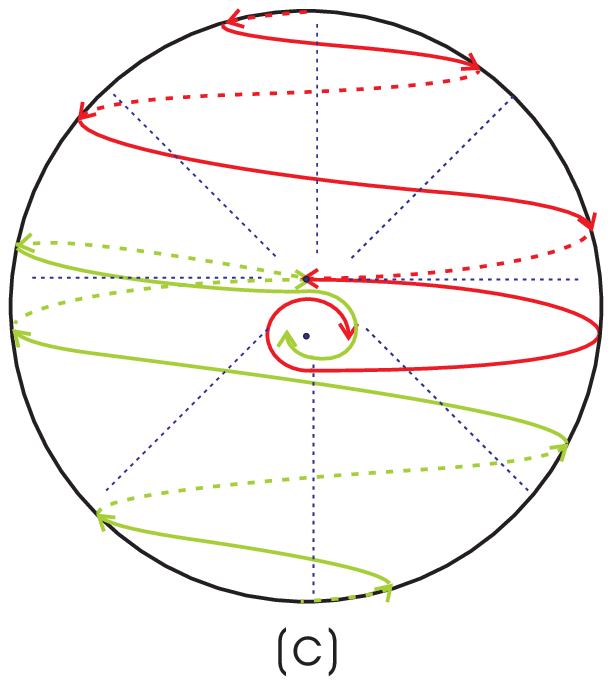} &
    \includegraphics[width=2.22in,height=2.2in]{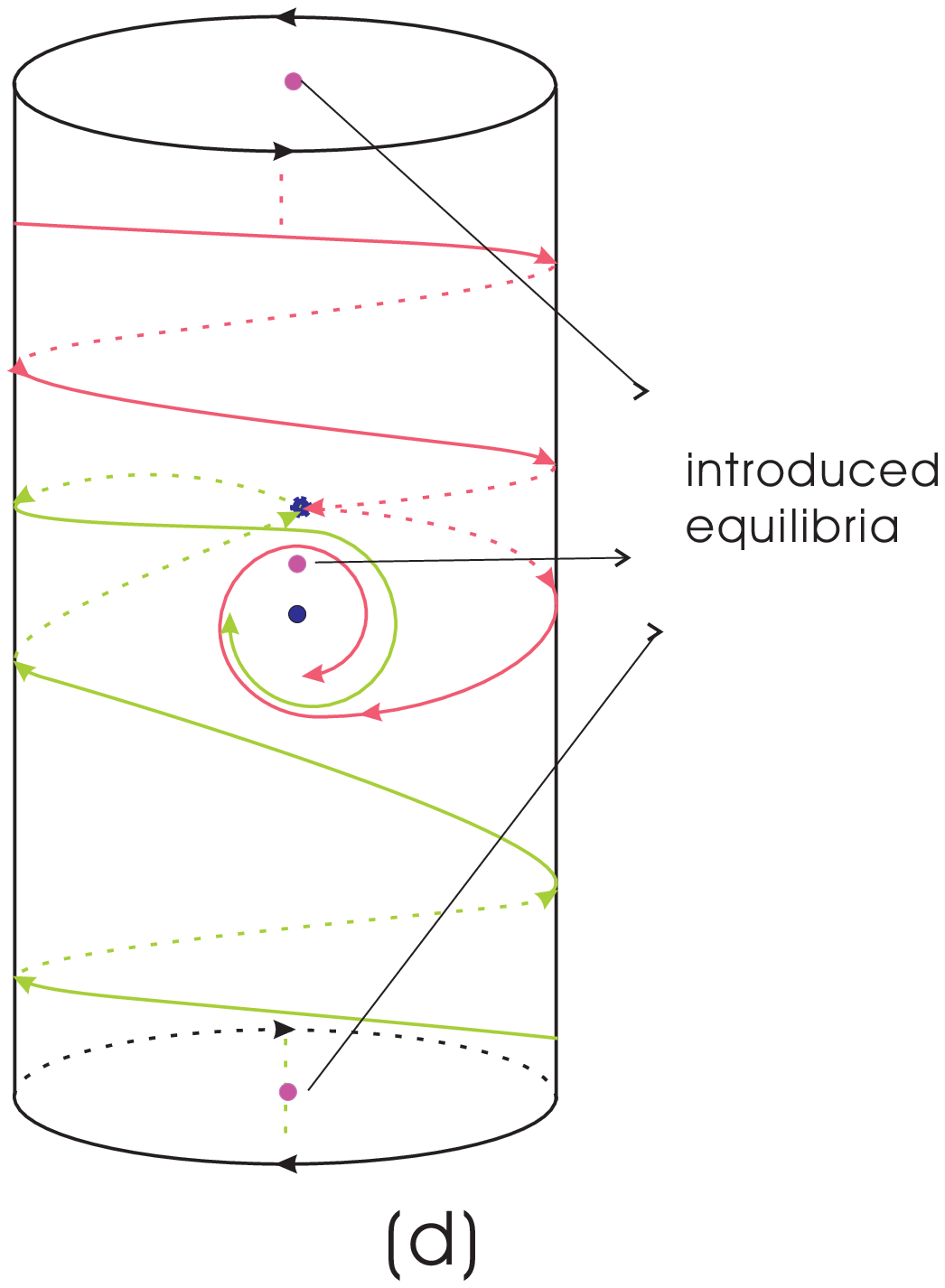}\\
    \end{tabular}
    \caption{Extending the dynamics through \emph{Klein bottle}}\label{klein bottle}
\end{center}
\end{figure}
\emph{By identifying $-\pi$ and $\pi$ and then gluing the two ends
together, we know that a pendulum is defined on a \emph{Klein
bottle} (see \textbf{Fig} \ref{klein bottle}.(b))}

\emph{Next open the nonorientable handle as stated in \textbf{Fig}
\ref{handle} from \textbf{Theorem}(\ref{extending dynamics}), the
surface dynamics can be effectively extended throughout the 3-ball
(\textbf{Fig} \ref{klein bottle}.(c)). Then after pulling and
expanding the two poles, (as shown in \textbf{Fig} \ref{klein
bottle}.(d)), we can glue the two ends back together and recover the
\emph{Klein bottle}. This time the system is situated on the whole
solid \emph{Klein bottle} with the surface dynamics stay unchanged.
}

\emph{From \textbf{Proposition}(\ref{identifying klein bottle}), we
know that there is exactly one nonorientable 3-manifold with a
\emph{genus} $1$ \emph{Heegaard Splitting}, and since \emph{Klein
bottle} is a nonorientable \emph{genus} $1$ surface, the
\emph{identity} map will certainly be the homeomorphism that glues
the two of them together. So in our pendulum case, there will be
exactly two same systems defined on the solid \emph{Klein bottle} in
the above way, and via the \emph{Heegaard diagram}, these two
3-manifolds will be glued by the \emph{identity} map obtaining a
nonorientable 3-manifold. }
\end{example}

\begin{example}
\emph{As shown in [Banks, 2002], a surface of \emph{genus} $2$ can
only carry two distinct knot types. \textbf{Fig} \ref{2 hole knot}
gives us the whole procedure of transforming a simplest knot to one
type of those which can be situated on a 2-hole surface by
performing the \emph{C-homeomorphisms}.}
\begin{figure}[!hbp]
  \begin{center}
  \begin{tabular}{cc}
  \includegraphics[width=2.22in,height=2.22in]{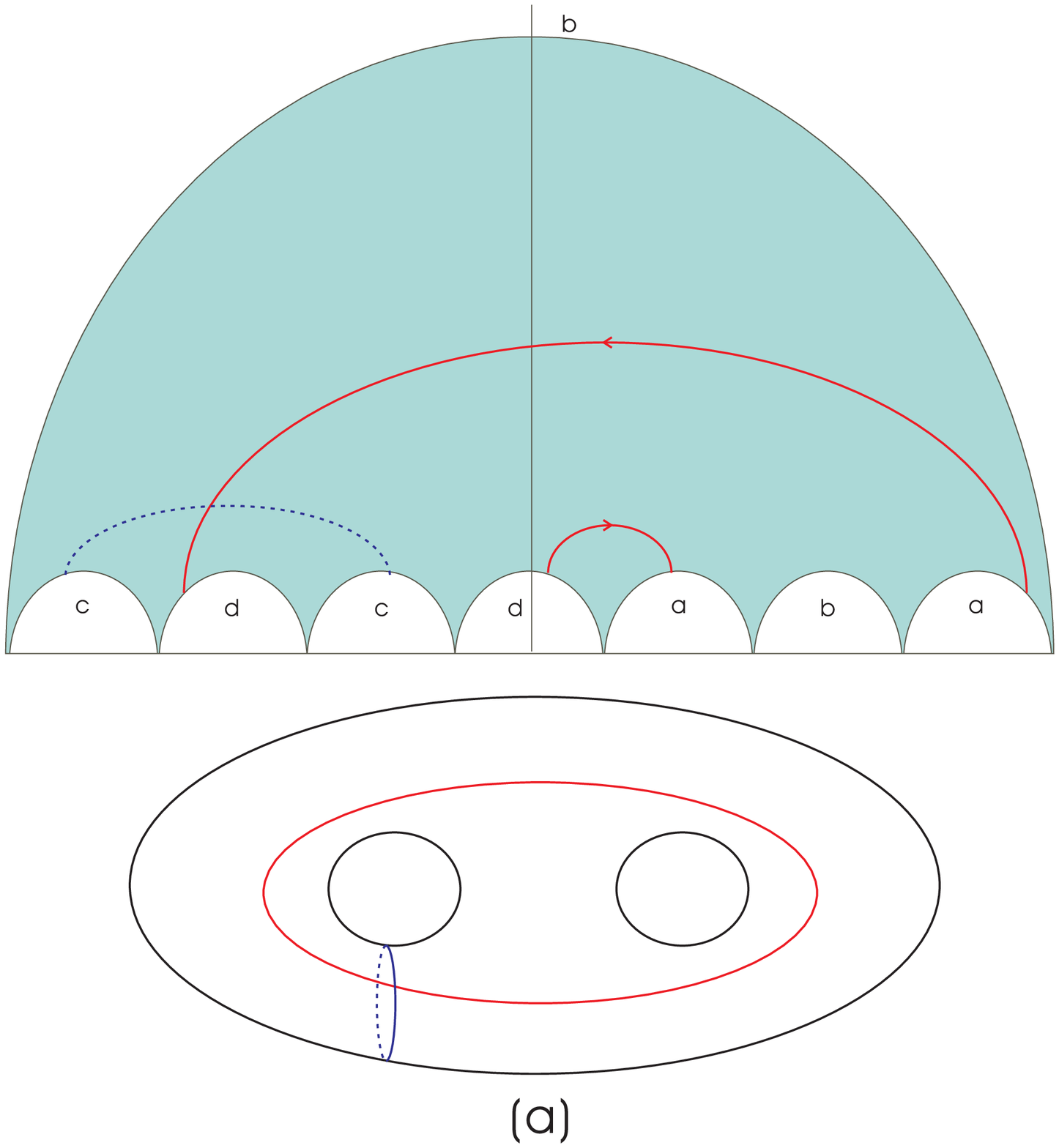} &
  \includegraphics[width=2.22in,height=2.22in]{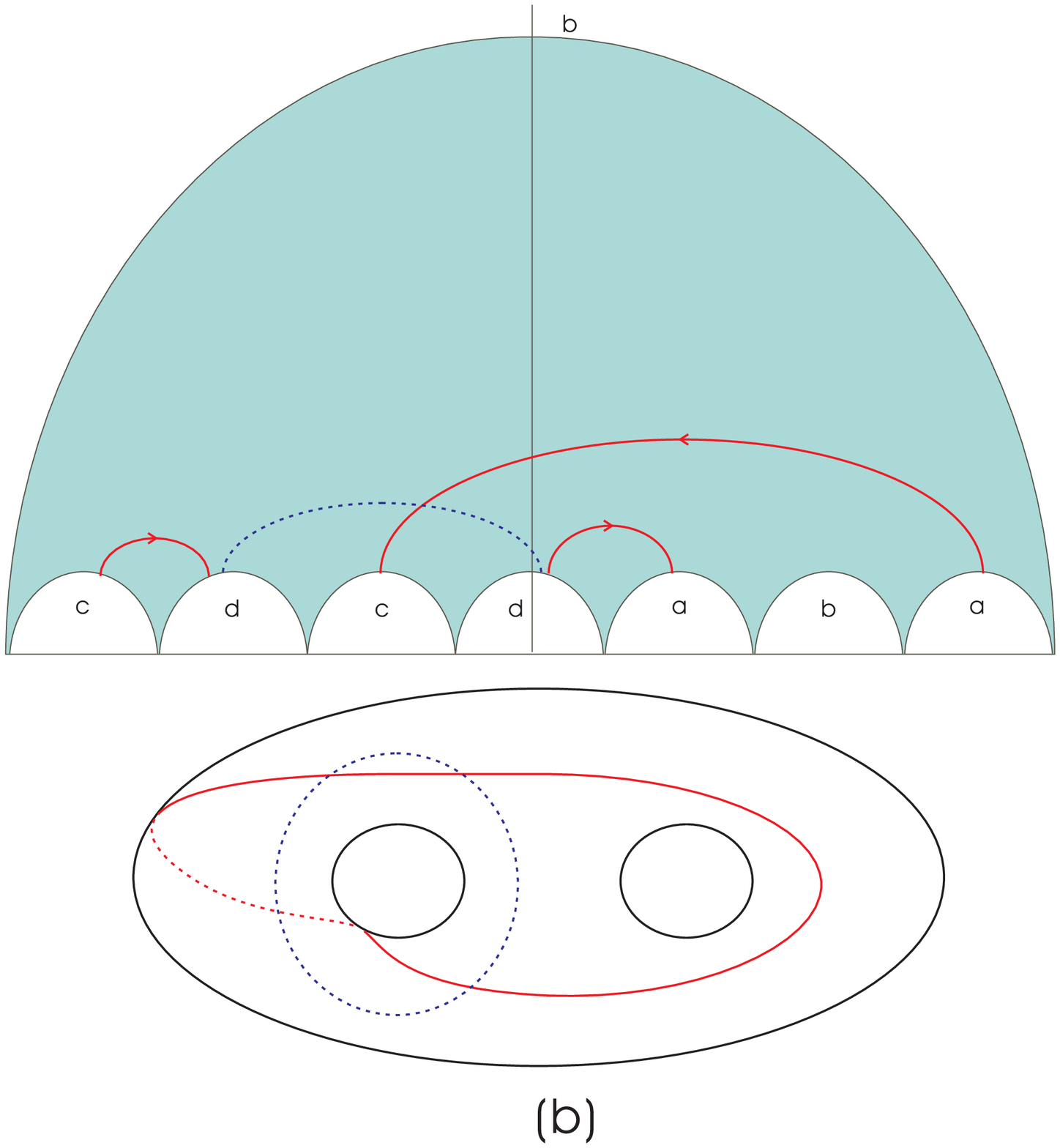}\\
  \includegraphics[width=2.22in,height=2.22in]{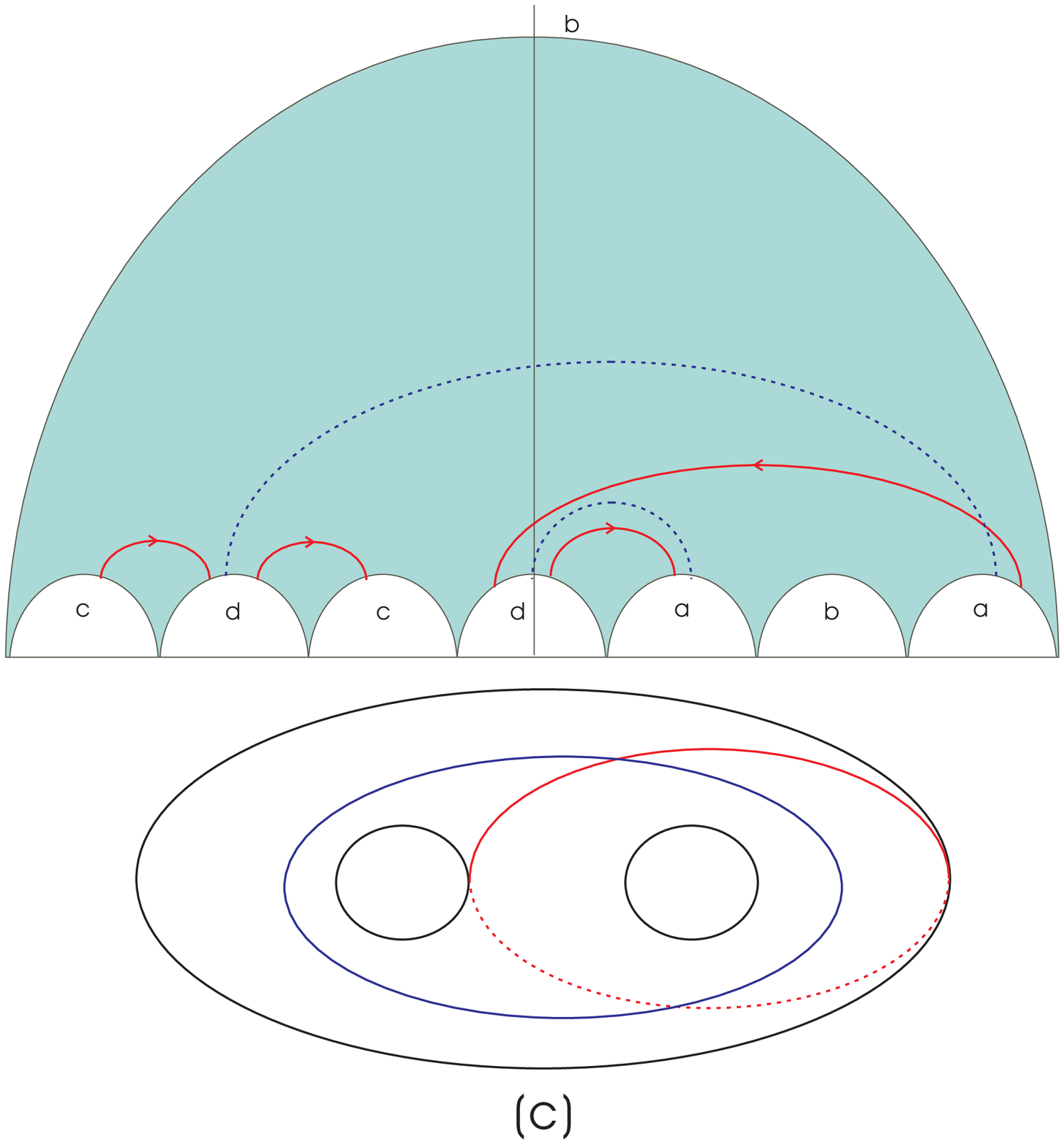} &
  \includegraphics[width=2.22in,height=2.22in]{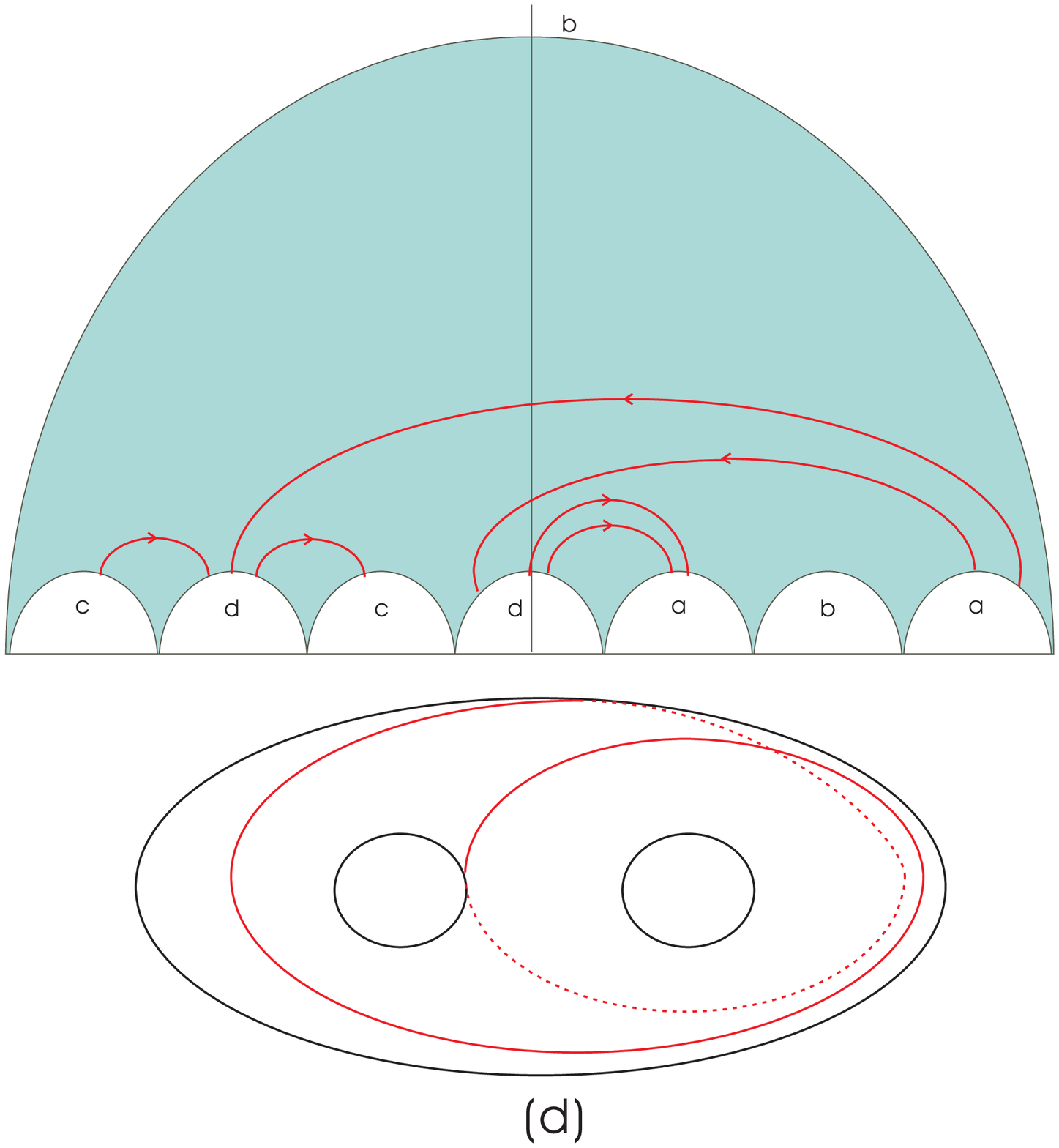}\\
  \end{tabular}
  \caption{Transforming a simple knot to one type of the two that
  can be carried by a surface of \emph{genus} $2$ via \emph{C-homeomorphisms}}
  \label{2 hole knot}
  \end{center}
\end{figure}

\emph{Also in [Banks \& Song, 2006], we gave an explicit
construction of a system that is situated on a \emph{$2$-hole torus}
by using the \emph{generalized automorphic functions}. The system
itself has \emph{Fuchsian group} generated by the transformations}
\begin{eqnarray*}
T_1(z)&=&-\frac{2z+13}{z+6} \\
T_2(z)&=&-\frac{1}{z+4} \\
T_3(z)&=&\frac{6z-13}{z-2} \\
T_4(z)&=&7z-28
\end{eqnarray*}

\emph{Choose}
\[H_1(z) = \frac{1}{z+2-3i},\quad H_2(z) = \frac{1}{z-2-3i}\]
\emph{we obtain a system with a \emph{pole} at $-2+3i$ and a
\emph{zero} at $2+3i$. The actual dynamics is shown in \textbf{Fig}
\ref{2 hole system}.(a)}
\begin{figure}[!hbp]
  \begin{center}
  \begin{tabular}{cc}
  \includegraphics[width=2.22in,height=1.8in]{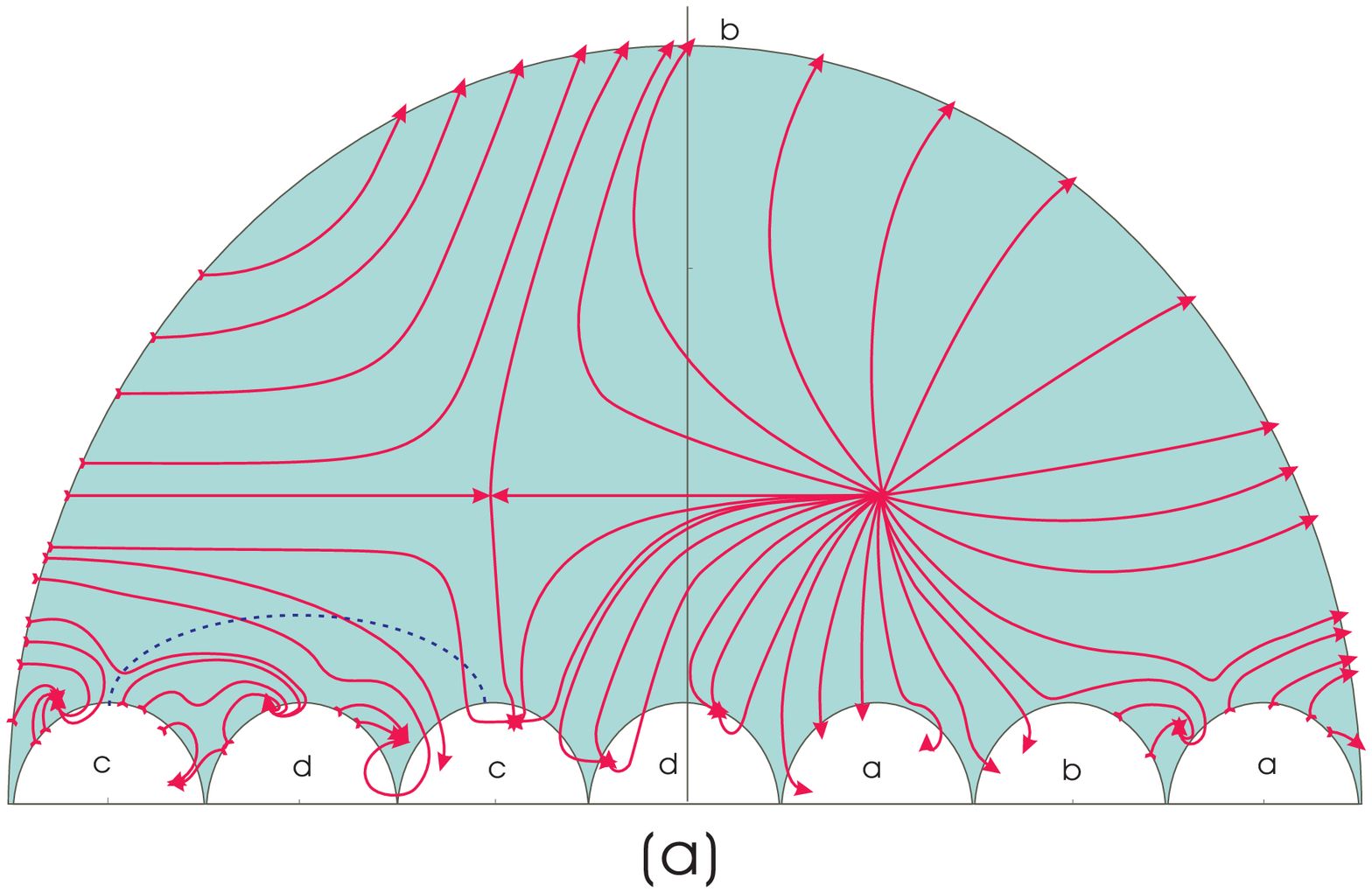} &
  \includegraphics[width=2.22in,height=1.8in]{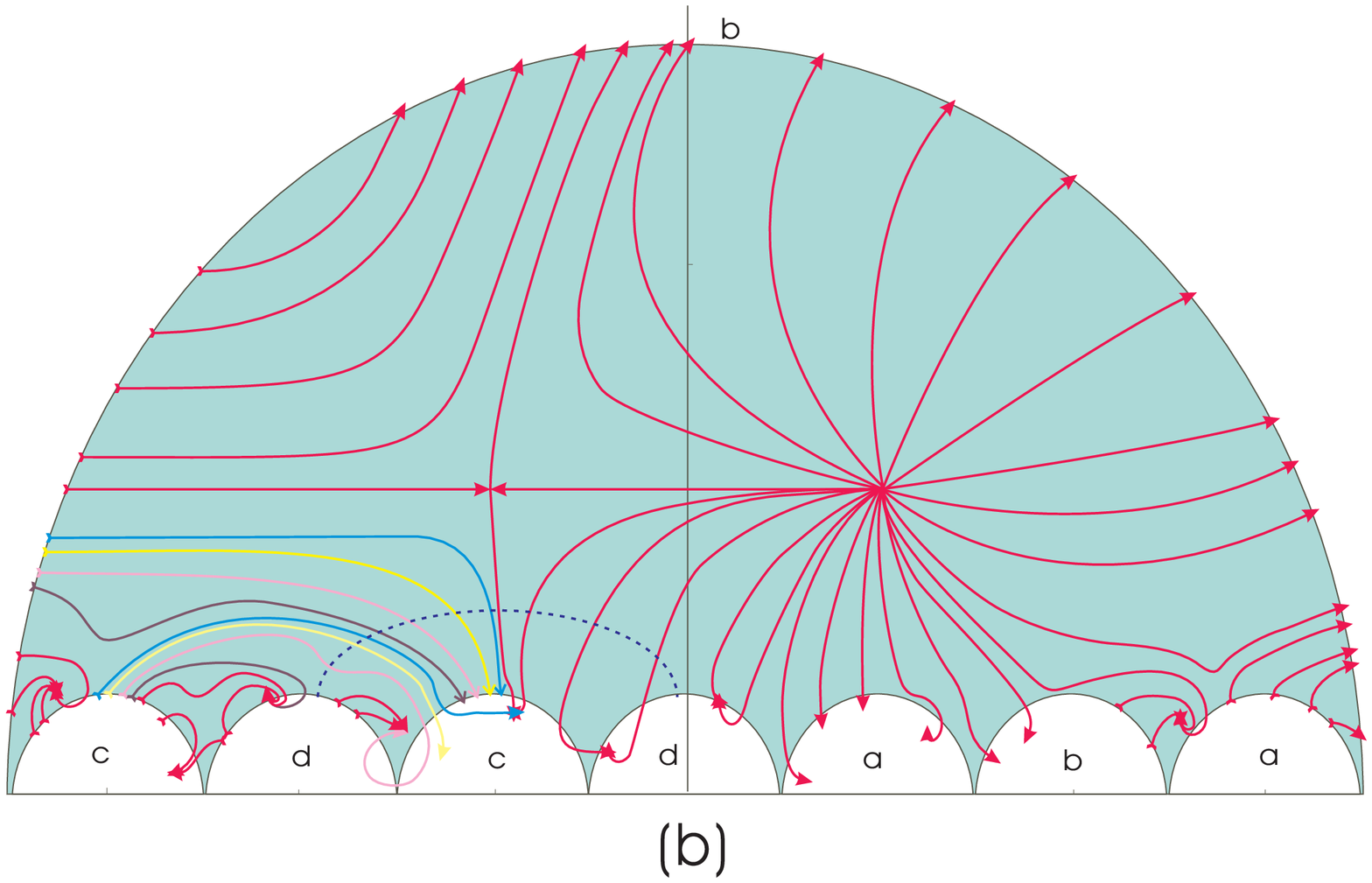}\\
  \includegraphics[width=2.22in,height=1.8in]{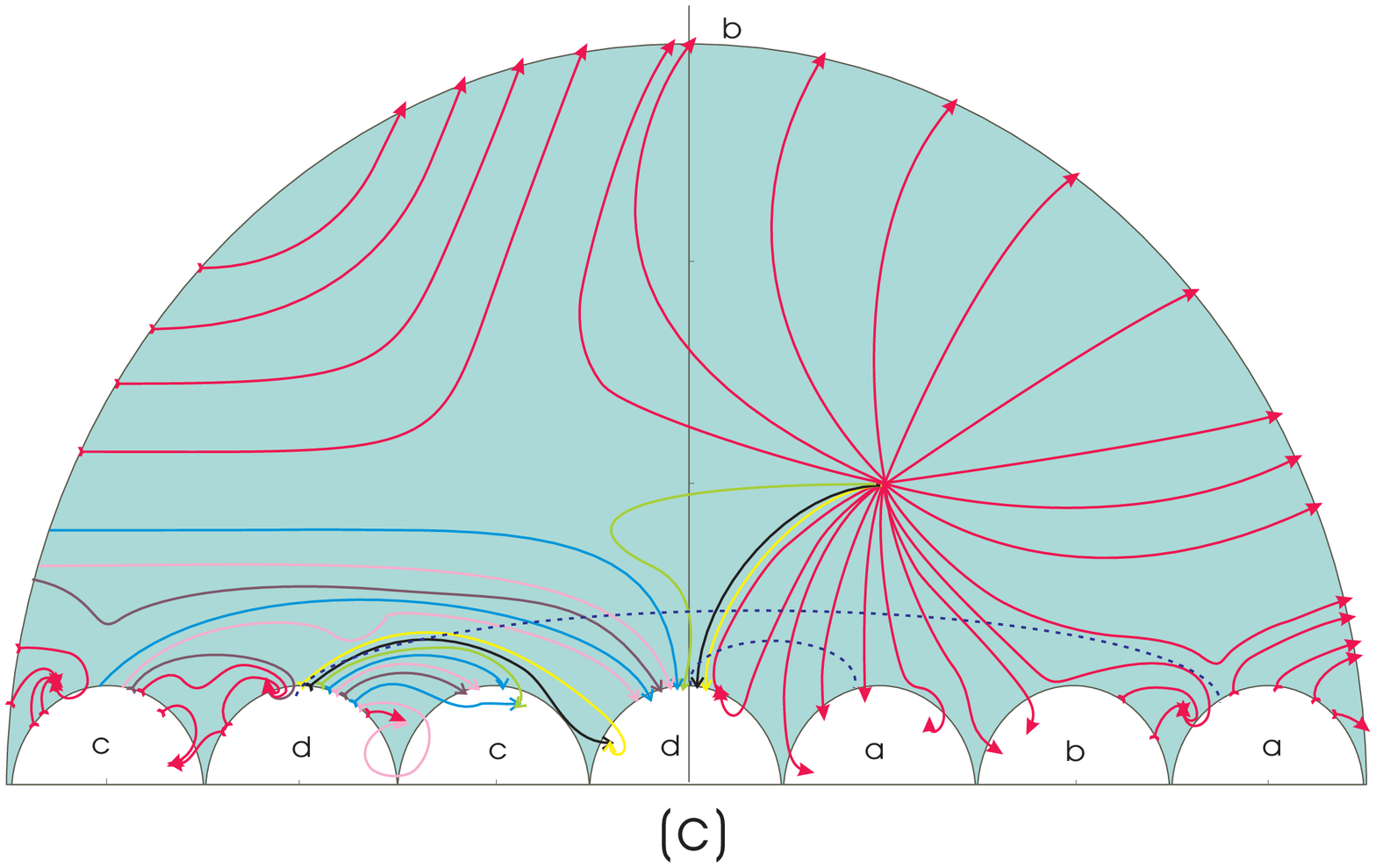} &
  \includegraphics[width=2.22in,height=1.8in]{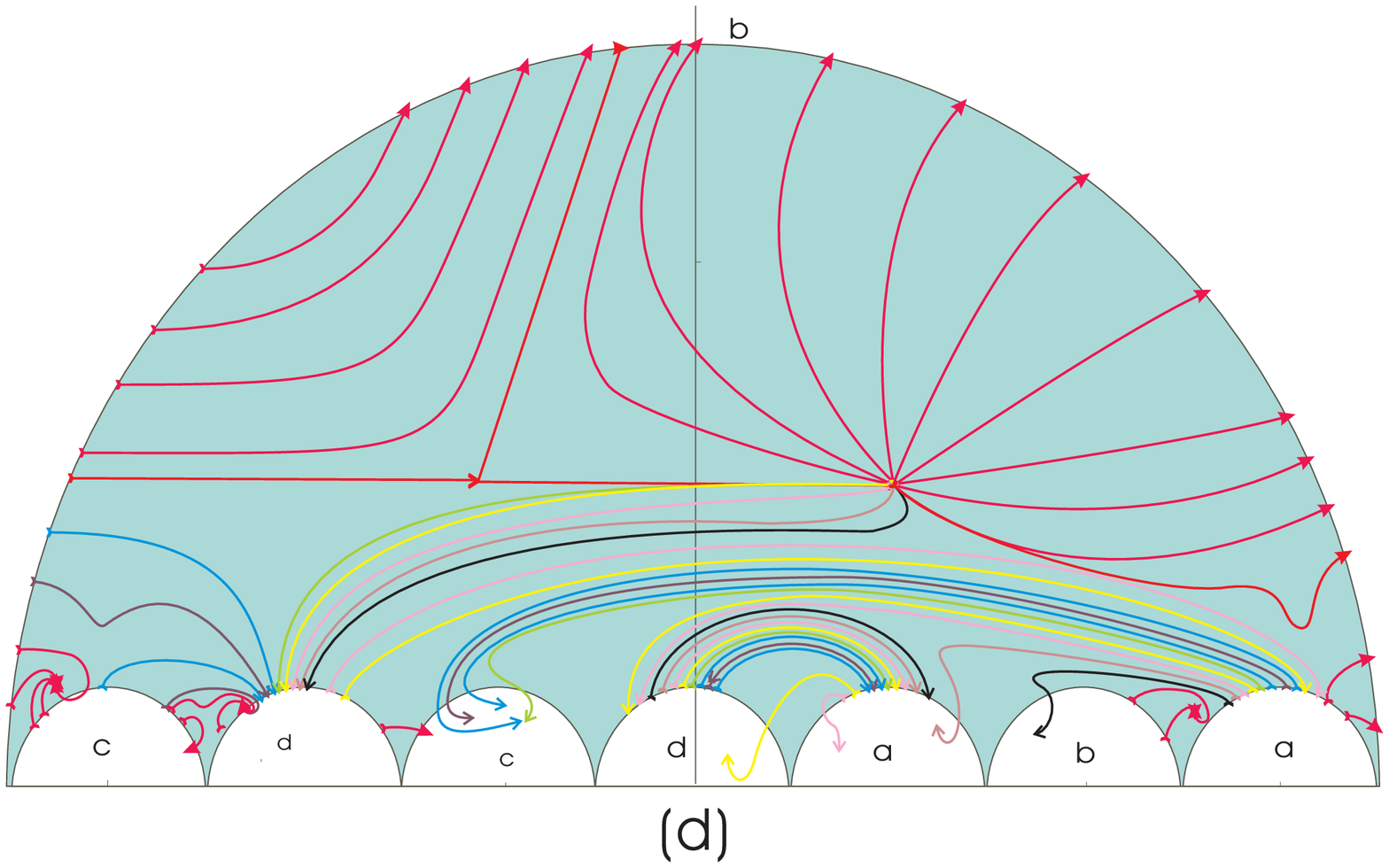}\\
  \end{tabular}
  \caption{Change of the surface dynamics according to the \emph{C-homeomorphism} surgery performed }
  \label{2 hole system}
  \end{center}
\end{figure}

\emph{When performing the surgery on the \emph{genus $2$} surface
(as shown in \textbf{Fig} \ref{2 hole knot}), the dynamics is also
changed accordingly. \textbf{Fig} \ref{2 hole system}.(b)-(d)
explain this procedure.}

\emph{We can then extend the system throughout the solid
$2$\emph{-hole torus} respectively, and use the
\emph{C-homeomorphisms} introduced in \textbf{Fig} \ref{2 hole
system} to glue the surface while the matching of the dynamics is
being guaranteed.}

\emph{In this way, we obtain a new system which is defined on a more
complicated 3-manifold from two simpler ones each sits on a solid
$2$\emph{-hole torus}.}

\end{example}

\section{Three-Dimensional Dynamical Systems and \emph{Heegaard
Splittings}}

In this section we consider a three-dimensional dynamical system
defined on a three-manifold without boundary containing only a
finite number of equilibria. we shall examine conditions under which
such a system has a \emph{Heegaard Splitting} that respects the
dynamics, i.e., contains an invariant \emph{genus p} surface, which
defines a \emph{Heegaard Splitting}. Our main technical tools will
be the \emph{Poincar$\acute{e}$-Hopf index} theorem and the
\emph{flow-box} theorem. The latter may be stated as follows:

\begin{theorem}
If $\phi_t$ denotes a dynamical system on a manifold M of dimension
n, then if $x \in M$ is not an equilibrium point (\emph{i.e.},
$\phi_t(x) \ne x, t \ne 0$), there exists a (closed) local
coordinate neighbourhood U of x such that on U, $\phi_{t}$ is
topological conjugate to the dynamical system
\begin{displaymath}
\left.\begin{array}{c}
\dot{x}_1=c\\
\dot{x}_{2}=0\\
\vdots\\
\dot{x}_{n}=0
\end{array} \right\}
\quad \quad x \in \{0 \le x_i \le 1, 1 \le i \le n \}
\end{displaymath}
where \emph{c} is a constant. \quad \qquad $\Box$
\end{theorem}

\begin{figure}[!hbp]
  \begin{center}
  \includegraphics[width=3in]{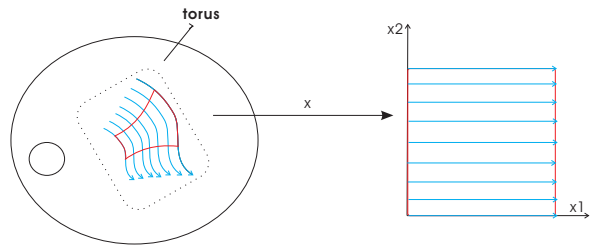}\\
  \caption{A local flowbox in 2-dimensional surface}\label{flowbox}
  \end{center}
\end{figure}

This says that locally, away from equilibria, the flow can be
``parallelized'', e.g., in two dimensions the flow looks locally
like the one in \textbf{Fig} \ref{flowbox}.

Since an invariant surface in \emph{M} can only have those
singularities of \emph{M}, in order that there exists an invariant
\emph{Heegaard Splitting} of \emph{genus} $p \ne 1$, the dynamical
system must have at least one equilibrium, so that systems with no
equilibria can only have \emph{genus $1$ Heegaard Splittings}, i.e.,
\emph{torus} or \emph{Klein bottle} splittings.

\begin{theorem}
In order that a 3-dimensional dynamical system on a compact manifold
M has a Heegaard splitting (compatible with the dynamics) of genus
p, it is necessary that it contains at least one equilibrium and
that in some subset of the equilibria, ${M_1, \cdots, M_l}$, there
is an invariant two-dimensional local surface passing through the
equilibrium with (2-dimensional) index $\vartheta_i$, such that

\[\sum_{i=1}^{l} \vartheta_i =2(1-p) \quad \quad \qquad \Box\]
\end{theorem}

\begin{corollary}
A dynamical system on a compact 3-manifold which has only
linearizable equilibria and a compatible Heegaard Splitting of genus
$p \ge 1$ must have at least $2(p-1)$ hyperbolic points. \quad \quad
\quad \quad \quad$\Box$
\end{corollary}

The above necessary conditions are not sufficient, in general, to
find sufficient conditions for a dynamical \emph{Heegaard Splitting}
we first recall the following result for a topological
\emph{Heegaard Splitting} and give a proof in order to motivate the
generalization.

\begin{theorem} \label{splitting}
(see [Hempel, 1976]) Every closed, connected 3-manifold M has a
Heegaard Splitting.
\end{theorem}
\textbf{Proof.} Take a triangulation \emph{K} of \emph{M} and let
$\Gamma_1$ be the set of all \emph{$1$-simplexes} of \emph{K}(i.e.,
the \emph{$1$-skeleton}). Let $\Gamma_2$ be the dual
\emph{$1$-skeleton}, which is the maximal \emph{$1$-subcomplex} of
the first derived complex $K^{\prime}$ which is disjoint from
$\Gamma_1$. Then if we put
\[V_i = N(\Gamma_i, K^{\prime \prime})\]
where \emph{N} is the normal neighbourhood of $\Gamma_i$ with
respect to $K^{\prime \prime}$ (the second derived of \emph{K}), it
can be shown that $(V_1,V_2)$ is a \emph{Heegaard Splitting} of
\emph{M}.
$\Box$ \\

It follows that any \emph{Heegaard Splitting} can be described in
this way. Suppose there is a dynamical \emph{Heegaard Splitting} of
a dynamical system on a closed connected manifold \emph{M}. Let
\emph{K} be a triangulation of \emph{M} determining the splitting as
in \textbf{Theorem}(\ref{splitting}). Then if $V_i = N( \Gamma_i,
K'')$ as above, $S = V_1 \cap V_2$ is a surface which is invariant
under the dynamics. Since \emph{M} is compact, we can cover \emph{M}
by a finite number of open sets $\{F_1,F_2, \cdots, F_L\}$ where
$F_i$ is a flow box if it does not contain an equilibrium point of
the dynamics or just a neighbourhood of such a point otherwise.
Suppose that $\{p_1, \cdots, p_k\}$ are equilibrium points of the
dynamics which belong to \emph{S}, and that $p_i \in F_i$ $(1 \le i
\le k)$. (This can always be done by renumbering the $F_i$'s.) Let
\[{E_i}^j = F_i \cap V_j \quad \quad 1 \le i \le k, 1 \le j \le 2\]
Then we can find a refinement $\{ {F_1}', \cdots, {F_{l_1}}',
{F_1}'', \cdots, {F_{l_2}}''\}$ of the remaining open sets
$\{F_{k+1}, \cdots, F_L\}$ so that there exists a partition
\[\Gamma^1 = \{E_1^1, E_2^1, \cdots, E_k^1, {F_1}', \cdots,{F_l}'\}\]
\[\Gamma^2 = \{E_1^2, E_2^2, \cdots, E_k^2, {F_1}'', \cdots,{F_l}''\}\]
such that
\[\cup \Gamma^i \subseteq V_i, \quad \quad 1 \le i \le 2\]
so that the sets $\Gamma^1$ and $\Gamma^2$ are invariant under the
dynamics. By taking the flow boxes small enough, we can associate a
triangulation of the manifold \emph{M} (by taking the corners of the
flow boxes away from the vertices) which is arbitrarily close to the
original one. Clearly, conversely, if we can find a system of flow
boxes for the dynamics on \emph{M} with the above properties and the
associated triangulation, then we will have a dynamical
\emph{Heegaard Splitting}. Thus we have proved
\begin{theorem}
Consider a compact 3-manifold M on which is given a compact
dynamical system. Suppose there is a refinement $\Gamma^1 \cup
\Gamma^2$ of a covering of M by flow boxes or neighbourhoods of
equilibria, such that $\Gamma^1$ and $\Gamma^2$ are invariant under
the dynamics. Let $\Gamma^1$ and $\Gamma^2$ be triangulations of
$\cup \Gamma^1, \cup \Gamma^2$, respectively, such that $\Gamma^1
\cup \Gamma^2$ is a triangulation of $M'$. Then $(\cup \Gamma^1,
\cup \Gamma^2)$ is a dynamical Heegaard Splitting of M if $\Gamma^1$
and $\Gamma^2$ are dual triangulations or the two-skeletons of
$\Gamma^1$ and $\Gamma^2$ have equal Euler characteristics.
\quad\quad $\Box$
\end{theorem}

\section{\emph{Connected Sums}}
\emph{Connected Sums} of 2- and 3-manifolds provide an effective
means of generating `complicated' manifolds out of simpler ones. In
this section we shall consider sums of dynamical systems on 2- and
3-manifolds.

Consider first the case of 2-manifolds. Given two (topological)
2-manifolds $S_1$ and $S_2$, their connected sum $S_1\#S_2$ is
obtained by removing discs $D_1$, $D_2$ from $S_1$ and $S_2$ and
sewing $S_1 / S_2$ to $S_2/ D_2$ along the boundaries of the discs.
If $S_2$ is a sphere, note that
\begin{equation}
S_1 \# S_2 = S_1
\end{equation}

\begin{figure}[!hbp]
\begin{center}
\includegraphics[width=3in]{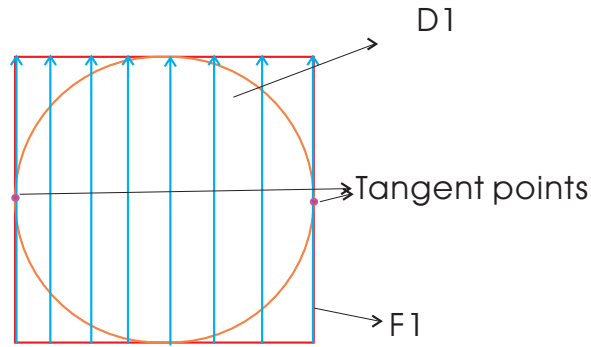}
\caption{Trajectories in flow box $F_1$ tangent to disc $D_1$}
\label{tangent}
\end{center}
\end{figure}

\begin{lemma}\label{nocriticalpoints}
Let $S_1$ and $S_2$ be two surfaces on which dynamical systems
$\phi_1$ and $\phi_2$ are defined. If we form the connected sum by
removing discs $D_1$, $D_2$ from $S_1$ and $S_2$ away from any
critical points, then we must introduce two hyperbolic equilibria
(with index $-1$) on the disc boundaries.
\end{lemma}
\textbf{Proof.} Since there are no equilibrium points in the discs
being removed, we can find flow boxes $F_1$ and $F_2$ in $S_1$ and
$S_2$, respectively, so that
\[D_i \le F_i, \quad \quad i=1, 2\]
provided $D_1, D_2$ are small enough. The discs can be chosen so
that there are two trajectories which are tangent to the discs at
two points. (see \textbf{Fig} \ref{tangent})

\begin{figure}[!hbp]
\begin{center}
\includegraphics[width=4in]{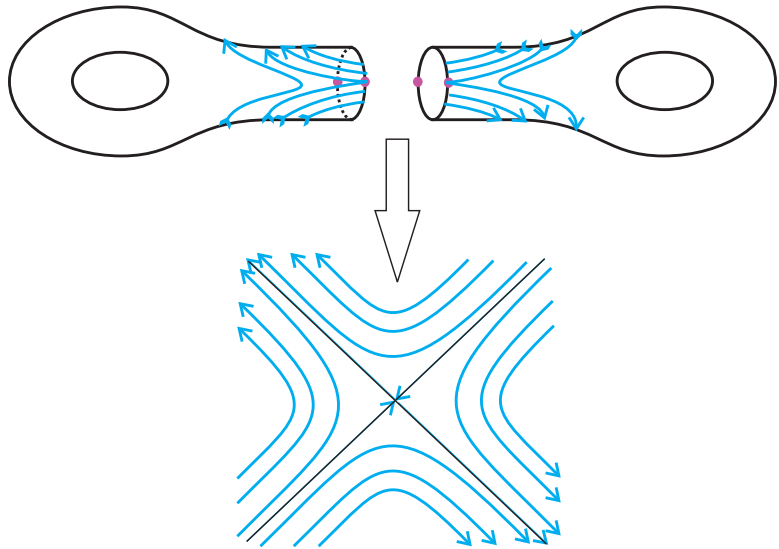}
\caption{Change to the dynamics after gluing two tori via
\emph{Connected Sum}} \label{connectedsum}
\end{center}
\end{figure}

If we now pull out tubes to form the \emph{connected sum} \emph{S},
these two points clearly become singular points on the
\emph{connected sum} as in \textbf{Fig} \ref{connectedsum}.

Suppose that one surface, say $S_2$, is a sphere. Since $S=S_1 \#
S_2=S_1$ in this case and $\chi (S_2) =2$, the total index of the
two singular points on the removed discs must be $-2$. \quad \quad
$\Box$

Suppose next that we form a \emph{connected sum} by removing discs
containing equilibria. There are several conditions that we have to
consider.

\begin{lemma} \label{withoutextra}
If we form the connected sum of two surfaces $S_1, S_2$ by removing
discs $D_1$ and $D_2$, which each contain an equilibrium $p_i
(i=1,2)$ without introducing new equilibria, then these equilibria
must be `dual' in the sense that if one equilibrium has $n_1$
elliptic sectors and $n_2$ hyperbolic sectors, then the other must
have $n_1$ hyperbolic sectors and $n_2$ elliptic ones.
\end{lemma}
\textbf{Proof.} Again we can assume $S_2$ is a sphere without loss
of generality. Let
\[S=S_1 \# S_2 =S_1\]
so that
\[ \chi(S) = \chi(S_1)\]
The index of one equilibrium point is
\[I(p_1) = 1+ \frac{n_1-n_2}{2}\]
Since $\chi(S_2) =2$, and without introducing extra equilibria,
$p_2$ must have index satisfying
\[I(p_2) +I(p_1) =2\]
so that
\[I(p_2) = 2- I(p_1) = 1+ \frac{n_2-n_1}{2}. \quad \quad \Box\]

Also, after removing discs containing equilibria and gluing the rest
together, we may introduce extra critical points on the disc
boundaries as well.

\begin{figure}
  \begin{center}
  \includegraphics[width=4in]{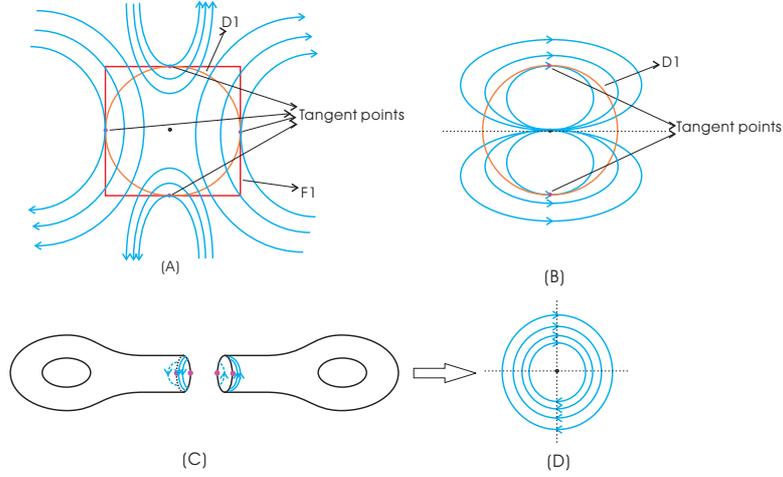}\\
  \caption{Gluing two tori while introducing extra equilibria}\label{saddleandcycle}
  \end{center}
\end{figure}

\begin{lemma}\label{samestructure}
Stick to the same notations as in
\textbf{Lemma}(\ref{withoutextra}), if there are new equilibria
being introduced, and the structure of $p_1$ and $p_2$ are exactly
the same, (i.e., $p_1$ and $p_2$ both have $n_1$ elliptic sectors
and $n_2$ hyperbolic sectors,) then the introduced equilibria must
be $n_1$ elliptic (with index $+1$) and $n_2$ hyperbolic (with index
$-1$).
\end{lemma}
\textbf{Proof.} Without loss of generality,we first look at the case
of a hyperbolic equilibrium (with index $-1$). It has $4$ hyperbolic
sectors. The removed discs $D_1$ and $D_2$ can be chosen such that
there are exactly four trajectories tangent to the discs at four
different points, as shown in \textbf{Fig} \ref{saddleandcycle}.(A)
.
The same argument in the proof of
\textbf{Lemma}(\ref{nocriticalpoints}) applies here. Referring to
\textbf{Fig} \ref{connectedsum}, one hyperbolic sector generates one
hyperbolic equilibrium (with index $-1$) after the gluing. And since
$p_i$ has $n_2$ hyperbolic sectors, we end up with $n_2$ hyperbolic
equilibria (all with index $-1$ respectively) being introduced after
the gluing via \emph{connected sum}.

We next consider the elliptic sectors. Suppose $p_1$ and $p_2$ only
contain $n_1$ elliptic sectors, as shown in \textbf{Fig}
\ref{saddleandcycle}.(B), within an elliptic sector, it is always
possible to find two closed elliptic trajectories which are tangent
to discs $D_1$ and $D_2$, respectively. If we then pull out tubes to
form the \emph{connected sum S}, these two points will certainly
turn into two singular points on the sphere, which are elliptic
equilibria and contain cycles only. (See \textbf{Fig}
\ref{saddleandcycle}.(C) and (D)).

Still, we let $S_2$ be a sphere and
\[S=S_1 \# S_2 =S_1\]
such that
\[ \chi(S) = \chi(S_1) = m\]
also we have
\[I(p_1) =I(p_2) =1+\frac{n_1}{2}\]
Suppose there are $n$ elliptic equilibria (with index $1$) being
introduced after the gluing, and since $\chi(S_2)=2$, we have
\[m-(1+\frac{n_1}{2})+2-(1+\frac{n_1}{2})+n= m\]
which gives us $n=n_1$. So there are $n_1$ new elliptic equilibria
appearing in $S$. \quad \quad \qquad $\Box$

Certainly, the structure of $p_1$ can be different from $p_2$ even
if there are extra equilibria being introduced.

\begin{lemma}
If we form the connected sum of two surfaces $S_1, S_2$ via removing
discs $D_1$ and $D_2$ which each contain an equilibrium $p_i
(i=1,2)$, then there must exist a separation to the sectors in
$p_1$: $n_{11}$ elliptic and $n_{21}$ hyperbolic sectors share the
same structure as those in $p_2$, while the rest are `dual' to the
remaining in $p_2$.
\end{lemma}
\textbf{Proof.} The proof follows from those of
\textbf{Lemma}(\ref{withoutextra}) and (\ref{samestructure}) since
these two are the only conditions that can happen to the dynamics
situated on surfaces when performing the \emph{connected sum}.
Separate $n_1$ elliptic sectors of $p_1$ to $n_{11}$ and $n_{12}$,
$n_2$ hyperbolic ones to $n_{21}$ and $n_{22}$, with $n_{11}$ and
$n_{21}$ being attached to the same structure on $p_2$, while
$n_{12}$ and $n_{22}$ being glued to their `dual' respectively.

Again, without loss of generality, we assume one surface, $S_2$, is
a sphere such that $\chi(S_2)=2$ and $S=S_1\#S_2=S_1=m$. And since
\begin{equation}\label{indexp1}
I(p_1) = 1+\frac{n_1-n_2}{2}=1+\frac{n_{11}+n_{12}-n_{21}-n_{22}}{2}
\end{equation}
\begin{equation}\label{indexp2}
I(p_2) =1+\frac{n_{11}+n_{22}-n_{21}-n_{12}}{2}
\end{equation}
From \textbf{Equation}(\ref{indexp1}), (\ref{indexp2}) and
\textbf{Lemma}(\ref{samestructure}),
\[m-I(p_1)+2-I(p_2)+n_{11} -n_{21}=m \]
is satisfied. \quad \quad \qquad $\Box$

We now extend the above results to the three-dimensional case. In
this case, the \emph{connected sum} of two compact 3-manifolds
$M_1,M_2$ is defined by removing two 3-cells from $M_1$ and $M_2$
and attaching their (spherical) boundaries together. This time, the
\emph{Euler Characteristic} of a compact 3-manifold is $0$, so by
$Poincar\acute{e}$-$Hopf$ theorem, the total index of any vector
field on the manifold is $zero$. First we form a \emph{connected
sum} by removing 3-cells which contain no equilibria. This time the
singular set is a (topological) circle, so we must introduce an
infinite set of equilibria or a limit cycle - we can do this by
twisting the cells before gluing. Note that the cycle does not
change the index, as expected. If we perform the connected sum by
removing cells containing equilibria without introducing new
singularities, then the equilibria must be `dual' in the sense that
regions on one part which point out of the cell must be matched by
those on the other part which point inwards. Clearly, the indices of
such critical points in 3-dimensions are the inverse of each other,
going a total index change of $0$, again as expected by the
$Poincar\acute{e}$-$Hopf$ theorem. If during the procedure of
removing 3-cells containing critical points, we introduce new
singularities, then from the combination of the statements above, we
know the total change of index is still $zero$.

\section{Conclusions}
In this paper, we have show how to generate a new dynamical system
on a complicated 3-manifold from a given one situated on a much
simpler 3-manifold by considering the corresponding \emph{Heegaard
diagram}, \emph{C-homeomorphisms} and the resulting dynamics on the
boundaries. Also, we gave the sufficient conditions under which a
system, which is defined on a three-manifold without boundary while
containing only a finite number of equilibria, has a \emph{Heegaard
Splitting} which respects the dynamics. A deeper look at
\emph{Connected Sum} and its effect on the natural dynamics will be
taken in the future paper.

\end{document}